# Practical Global Optimization Algorithm for the Sum-of-Ratios Problem

**Yunchol Jong** [a]

[a] Center of Natural Sciences, University of Science, Pyongyang, DPR Korea

**July 22, 2012**

**Abstract.** This paper presents a practical method for finding the globally optimal solution to the sum-of-ratios problem arising in image processing, engineering and management. Unlike traditional methods which may get trapped in local minima due to the non-convex nature of this problem, our approach provides a theoretical guarantee of global optimality. Our algorithm is based on solving a sequence of convex programming problems and has global linear and local superlinear/quadratic rate of convergence. The practical efficiency of the algorithm is demonstrated by numerical experiments for synthetic data.

## 1. Introduction

The sum-of-ratios problem, which is to minimize a sum of several fractional functions subjected to convex constraints, is a non-convex optimization problem that is difficult to solve by traditional optimization methods. The problem arises in many applications such as optimization of the average element shape quality in the finite element method, computer graphics and management ([1]~[3]). In [3], many problems of projective geometry including multiview triangulation, camera resectioning and homography estimation have been formulated as the sum-of-ratios problem and a branch-and-bound method has been proposed to find its global solution which relies on recent developments in fractional programming and the theory of convex underestimators ([5]~[7]). In the method of [3], number of variables increases as twice as the number of fractional functions involving in the sum and a second-order cone programming problem is needed to solve for obtaining a lower bound of the optimal value in each iteration. Their algorithm is provably optimal, that is, given any tolerance $\varepsilon$, if the optimization problem is feasible, the algorithm returns a solution which is at most $\varepsilon$ far from the global optimum. The branch-and-bound method requires a lot of computations, has low convergence and it is not easy to find a reasonable branching strategy.

Recently there has been some progress made towards finding the global solution to a few of these optimization problems ([8], [9]). However, the resulting algorithm is numerically unstable, computationally expensive and does not generalize for more views or harder problems like resectioning. In [10], linear matrix inequalities were used to approximate the global optimum, but no guarantee of actually obtaining the global optimum is given. Also, there are unsolved problems concerning numerical stability. Robustification using the $L_1$-norm was presented in [11], but the approach is restricted to the affine camera model. In [12], a wider class of geometric reconstruction problems was solved globally, but with $L_1$-norm.

In this paper, an efficient algorithm is presented which transforms the sum-of-ratios problem into parametric convex programming problem and finds the global solution successfully.

E-mail: yuncholjong@yahoo.com

## 2. Transformation into the equivalent parametric convex programming

The sum-of ratios problem seeks to minimize the sum of $N$ fractional functions subject to convex constraints, which is formulated as follows.

$$\min F(x) = \sum_{i=1}^{N} F_i(x),$$
$$\text{subject to } g_i(x) \leq 0, \ i = 1,...,m, \quad (1)$$
$$x \in R^n$$

where $F_i(x) = \dfrac{f_i(x)}{h_i(x)}$, $i = 1,...,N$, and $f_i(x)$, $g_i(x)$, $-h_i(x)$ are continuously differentiable convex functions, respectively. It is assumed that $f_i(x) \geq 0$ and $h_i(x) > 0$ in the feasible set $X = \{x \in R^n \mid g_i(x) \leq 0, i = 1,...,m\}$ and $IntX = \{x \in R^n \mid g_i(x) < 0, 1 \leq i \leq m\} \neq \emptyset$. Even with these restrictions the above problem is NP-complete [4].

It is easy to see that the problem (1) is equivalent to the following problem.

$$\min \sum_{i=1}^{N} \beta_i,$$
$$\text{subject to } F_i(x) \leq \beta_i, \ i = 1,...,N \quad (2)$$
$$g_i(x) \leq 0, \ i = 1,...,m$$
$$x \in R^n$$

Let us consider the optimality condition of this problem. The Lagrange function of the problem (2) is $L(x, \beta, u, v) = \sum_{i=1}^{N} \beta_i + \sum_{i=1}^{N} u_i (f_i(x) - \beta_i h_i(x)) + \sum_{i=1}^{m} v_i g_i(x)$ and the Karush-Kuhn-Tucker (KKT) condition for the problem (2) is as follows.

$$\frac{\partial L}{\partial x} = \sum_{i=1}^{N} u_i (\nabla f_i(x) - \beta_i \nabla h_i(x)) + \sum_{i=1}^{m} v_i \nabla g_i(x) = 0 \quad (3)$$

$$\frac{\partial L}{\partial \beta_i} = 1 - u_i h_i(x) = 0, \ i = 1,...,N \quad (4)$$

$$u_i \frac{\partial L}{\partial u_i} = u_i (f_i(x) - \beta_i h_i(x)) = 0 \quad (5)$$

$$v_i \frac{\partial L}{\partial v_i} = v_i g_i(x) = 0 \quad (6)$$

$$g_i(x) \leq 0, \ u_i \geq 0, v_i \geq 0, i = 1,...,m \quad (7)$$

Since $h_i(x) > 0$, $i = 1,...,N$ for every $x \in X$, (4) is equivalent to

$$u_i = \frac{1}{h_i(x)}, \ i = 1,...,N, \quad (8)$$

and so (5) is equivalent to

$$f_i(x) - \beta_i h_i(x) = 0, \quad i = 1, \ldots, N \tag{9}$$

The system (3), (6) and (7) is just the KKT condition of the problem for parameters $u_i$ and $\beta_i$, $i = 1, \ldots, N$.

$$\begin{aligned}
\text{Min} \quad & \sum_{i=1}^{N} u_i \left( f_i(x) - \beta_i h_i(x) \right), \\
\text{subject to} \quad & g_i(x) \leq 0, \quad i = 1, \ldots, m \\
& x \in R^n
\end{aligned} \tag{10}$$

For fixed $u_i > 0$, $\beta_i \geq 0, i = 1, \ldots, N$, the objective function of the problem (10) is convex and the problem (10) is a convex programming problem.
Therefore, we conclude that the solution of the problem (1) can be obtained by finding those satisfying (8), (9) among the solutions of the convex programming problem (10). If such solution is unique, the solution is just the global solution of the problem (1).
Let $\alpha = (\beta, u)$ denote parameter vector and let

$$\Omega = \{ \alpha = (\beta, u) \in R^{2N} \mid \beta \geq 0, u \geq l > 0 \}.$$

Let $x(\alpha)$ be the solution of the problem (10) for fixed $\alpha \in \Omega$ and let

$$\varphi(\alpha) = \sum_{i=1}^{N} u_i \left( f_i(x(\alpha)) - \beta_i h_i(x(\alpha)) \right) \tag{11}$$

$$\varphi_i(\alpha) = u_i \left( f_i(x(\alpha)) - \beta_i h_i(x(\alpha)) \right) \tag{12}$$

$$\bar{\psi}_i^1(\alpha) = f_i(x(\alpha)) - \beta_i h_i(x(\alpha)) \tag{13}$$

$$\bar{\psi}_i^2(\alpha) = 1 - u_i h_i(x(\alpha)) \tag{14}$$

If there exists a unique solution of the problem (10) for each fixed $\alpha \in \Omega$, then functions (11)-(14) are differentiable with $\alpha$ and we have

$$\frac{\partial \bar{\psi}_i^1(\alpha)}{\partial \beta_j} = \begin{cases} -h_i(x(\alpha)), & i = j \\ 0, & i \neq j \end{cases} \tag{15}$$

$$\frac{\partial \bar{\psi}_i^2(\alpha)}{\partial u_j} = \begin{cases} -h_i(x(\alpha)), & i = j \\ 0, & i \neq j \end{cases} \tag{16}$$

Then our problem (1) is reduced to find the solution of the $x(\alpha)$ which is also the solution of the system of nonlinear equations

$$\psi_i^1(\alpha) = -f_i(x(\alpha)) + \beta_i h_i(x(\alpha)) = 0, \quad i = 1, \ldots, N \tag{17}$$

$$\psi_i^2(\alpha) = -1 + u_i h_i(x(\alpha)) = 0, \quad i = 1, \ldots, N \tag{18}$$

## 3. Algorithm and its convergence.

Let
$$\psi_i(\alpha) = \psi_i^1(\alpha), i = 1,...,N$$
$$\psi_{N+i}(\alpha) = \psi_i^2(\alpha), i = 1,...,N$$

Then (17), (18) can be rewritten as

$$\psi(\alpha) = 0 \qquad (19)$$

**Lemma 1.** If $\psi(\alpha)$ is differentiable, it is strongly monotone with constant $\delta>0$ in $\Omega$, where
$$\delta = \min_i \delta_i, \quad \delta_i = \min_{x \in X} h_i(x)$$

**(Proof)** By (15) and (16), the Jacobian matrix of $\psi(\alpha)$ is as follows.

$$\psi'(\alpha) = \begin{bmatrix} h_1(x(\alpha)) & 0 & \cdots & 0 & 0 & 0 & \cdots & 0 \\ 0 & h_2(x(\alpha)) & \cdots & 0 & 0 & 0 & \cdots & 0 \\ \vdots & \vdots & \ddots & \vdots & \vdots & \vdots & \cdots & \vdots \\ 0 & 0 & \cdots & h_N(x(\alpha)) & 0 & 0 & \cdots & 0 \\ 0 & 0 & \cdots & 0 & h_1(x(\alpha)) & 0 & \cdots & 0 \\ 0 & 0 & 0 & 0 & 0 & h_2(x(\alpha)) & \cdots & 0 \\ \vdots & \vdots & \cdots & \vdots & \vdots & \vdots & \ddots & \vdots \\ 0 & 0 & 0 & 0 & 0 & 0 & \cdots & h_N(x(\alpha)) \end{bmatrix} \qquad (20)$$

Since $x(\alpha) \in X$ and $h_i(x(\alpha)) > 0$, $i=1,...,N$, $\psi'(\alpha)$ is positive definite. Therefore, for any $d \in R^{2N}$ we have

$$d^T \psi'(\alpha) d = \sum_{i=1}^N d_i^2 h_i(x(\alpha)) + \sum_{i=1}^N d_{i+N}^2 h_i(x(\alpha)) =$$
$$= \sum_{i=1}^N (d_i^2 + d_{i+N}^2) h_i(x(\alpha)) \geq \sum_{i=1}^N (d_i^2 + d_{i+N}^2) \delta_i \geq$$
$$\geq \delta \sum_{i=1}^{2N} d_i^2 = \delta \|d\|^2,$$

which completes the proof. □

Let
$$A_i(\alpha) = \frac{f_i(x(\alpha))}{h_i(x(\alpha))}, A_{i+N}(\alpha) = \frac{1}{h_i(x(\alpha))}, \quad i=1,...,N$$

**Lemma 2.** The equation $\psi(\alpha)=0$ is equivalent to the equation $\alpha = A(\alpha)$. Moreover, if the problem (1) has an optimal solution, the equation (19) has at least one solution in $\Omega$.

**(Proof)** The first proposition is obvious from the definition of $A(\alpha)$. As shown in section 2, if the problem (1) has an optimal solution, there is $\alpha \in \Omega$ such that the solution is also the solution of the problem (10) satisfying (19) for the $\alpha \in \Omega$. Thus, there is such $\alpha \in \Omega$ that $\alpha = A(\alpha)$. □

Let $B(\alpha) = \pi_\Omega(\alpha - \lambda \psi(\alpha))$, where $\pi_\Omega(a)$ is the projection operator of $a$ onto $\Omega$.

**Lemma 3.** Assume that $\psi(\alpha)$ is differentiable and satisfies the Lipschitz condition in $\Omega$. Then $B: \Omega \to \Omega$ is a contractive mapping.

**(Proof)** $\psi(\alpha)$ is strongly monotone by Lemma 1 and then it proves the proposition of the Lemma together with Lipschitz property and the non-expansivity of the projection. □

**Theorem 1.** Assume that $\psi(\alpha)$ is differentiable and satisfies the Lipschitz condition in $\Omega$. The equation $B(\alpha) = \alpha$ is equivalent to the equation $\psi(\alpha) = 0$ and the equation (19) has a unique solution.

**(Proof)** By Lemma 3 and the contractive mapping principle, $B(\alpha)$ has only one fixed point $\alpha^*$ in $\Omega$. From the property of projection, we have

$$\left[\pi_\Omega(\alpha^* - \lambda \psi(\alpha^*)) - (\alpha^* - \lambda \psi(\alpha^*))\right]^T \left[\alpha - \pi_\Omega(\alpha^* - \lambda \psi(\alpha^*))\right]$$
$$= \left[\alpha^* - (\alpha^* - \lambda \psi(\alpha^*))\right]^T (\alpha - \alpha^*) \geq 0$$

for $\alpha \in \Omega$, i.e. for every $\alpha \in \Omega$, we have

$$\psi(\alpha^*)^T (\alpha - \alpha^*) \geq 0. \tag{21}$$

On the other hand, there exists $\bar{\alpha} \in \Lambda \subset \Omega$ such that $\psi(\bar{\alpha}) = 0$ by Lemma 2 and so it follows from (21) that

$$(\psi(\alpha^*) - \psi(\bar{\alpha}))^T (\bar{\alpha} - \alpha^*) \geq 0 \tag{22}$$

By Lemma 1, $\psi(\alpha)$ is strongly monotone and so we have

$$(\psi(\alpha^*) - \psi(\bar{\alpha}))^T (\alpha^* - \bar{\alpha}) \geq \delta \|\alpha^* - \bar{\alpha}\|^2$$

This inequality and (22) together implies

$$0 \leq (\psi(\alpha^*) - \psi(\bar{\alpha}))^T (\bar{\alpha} - \alpha^*) \leq -\delta \|\alpha^* - \bar{\alpha}\|^2$$

Therefore, $\|\alpha^* - \bar{\alpha}\| = 0$, i.e. $\bar{\alpha} = \alpha^*$. Thus $\psi(\alpha^*) = 0$, $\alpha^* \in \Lambda$, which means that $\alpha^*$ is a solution of (19). Then $\alpha^*$ is a fixed point of $B(\alpha)$ in $\Omega$ and it is unique. Hence, $\alpha^*$ is a unique solution of the equation (19). □

**Theorem 2.** If $\psi(\alpha)$ is differentiable and satisfies the Lipschitz condition in $\Omega$, then the projection method

$$\alpha^{k+1} = \pi_\Omega(\alpha^k - \lambda \psi(\alpha^k))$$

converges to a solution of the equation (19) with linear rate.

**(Proof)** By Lemma 3, $\psi(\alpha)$ is contractive and then for the fixed point $\alpha^*$ of $B(\alpha)$ there is a $q \in (0,1)$ such that

$$\left\|\alpha^{k+1} - \alpha^*\right\| = \left\|B(\alpha^k) - B(\alpha^*)\right\| \leq q\left\|\alpha^k - \alpha^*\right\|,$$

which proves the linear convergence of $\{\alpha_k\}$ to $\alpha^*$. □

As seen in Lemma 2, $\alpha$ is a fixed point of $A(\alpha)$ if and only if $\alpha$ is a root of $\psi(\alpha)$. We can construct the following simple iterative method to find a fixed point of $A(\alpha)$.

$$\beta^{k+1} = \left(\frac{f_1(x^k)}{h_1(x^k)}, \cdots, \frac{f_N(x^k)}{h_N(x^k)}\right)^T, \quad u^{k+1} = \left(\frac{1}{h_1(x^k)}, \cdots, \frac{1}{h_N(x^k)}\right)^T \quad (23)$$

where $x^k = x(\alpha^k)$.

The Newton method for the equation (19) is as following.

$$\alpha^{k+1} = \alpha^k - \left[\psi'(\alpha^k)\right]^{-1}\psi(\alpha^k) \quad (24)$$

By (20), (17) and (18), the right-hand side of (24) is equal to $A(\alpha^k)$, i.e. the right-hand side of (23). Therefore, (24) means $\alpha^{k+1} = A(\alpha^k)$, that is, the simple iterative method to find a fixed point of $A(\alpha)$ is just the Newton method to solve the equation (19). Hence, the algorithm (23) has local superlinear or quadratic convergence rate.

**Theorem 3.** Assume that $\psi(\alpha)$ is differentiable, satisfies the Lipschitz condition and $\psi'(\alpha)$ satisfies the Lipschitz condition in $\Omega$, i.e.

$$\exists L, \forall \alpha, \alpha' \in \Omega, \ \left\|\psi'(\alpha) - \psi'(\alpha')\right\| \leq L\|\alpha - \alpha'\| \quad (25)$$

and that

$$\exists m, \ \forall \alpha \in \Omega, \ \left\|\left[\psi'(\alpha)\right]^{-1}\right\| \leq m. \quad (26)$$

Let

$$\alpha^{k+1} = \alpha^k + \lambda_k p^k, \ p^k = -\left[\psi'(\alpha^k)\right]^{-1}\psi(\alpha^k), \quad (27)$$

where $\lambda_k$ is the greatest $\xi^i$ satisfying

$$\left\|\psi(\alpha^k + \xi^i p^k)\right\| \leq (1 - \varepsilon\xi^i)\left\|\psi(\alpha^k)\right\| \quad (28)$$

and $i \in \{0,1,2,...\}$, $\xi \in (0,1)$, $\varepsilon \in (0,1)$.

Then, the modified Newton method defined by (27) and (28) converges to the unique solution $\alpha^*$ of $\psi(\alpha) = 0$ with linear rate for any starting point $\alpha^0 \in \Omega$ and the rate in the neighborhood of the solution is quadratic.

**(Proof)** We have already shown the existence and uniqueness of the solution to the equation $\psi(\alpha) = 0, \alpha \in \Omega$ above.

If there is $k$ such that $\psi(\alpha^k) = 0$, then $\alpha^k$ is a solution. So, it is assumed that $\psi(\alpha^k) \neq 0$ for every $k$.

For $\lambda \in [0,1]$, we have the following by the Newton-Leibnitz formula and (25).

$$\left\| \psi(\alpha^k + \lambda p^k) \right\| = \left\| \psi(\alpha^k) + \lambda \int_0^1 \psi'(\alpha^k + \theta \lambda p^k) p^k d\theta \right\| =$$
$$= \left\| \psi(\alpha^k) + \lambda \int_0^1 [\psi'(\alpha^k + \theta \lambda p^k) - \psi'(\alpha^k)] p^k d\theta - \lambda \psi(\alpha^k) \right\| \quad (29)$$
$$\leq (1-\lambda) \left\| \psi(\alpha^k) \right\| + \lambda^2 L \left\| p^k \right\|^2,$$

where we took the advantage of relation

$$\left\| \psi'(\alpha^k + \theta \lambda p^k) - \psi'(\alpha^k) \right\| \leq L \theta \lambda \left\| p^k \right\|$$

by the Lipschitz condition, and $\psi'(\alpha^k) p^k = -\psi(\alpha^k)$ by (27).

In view of (26), it follows from (27) and (29) that

$$\left\| \psi(\alpha^k + \lambda p^k) \right\| \leq \left[ 1 - \lambda \left( 1 - \lambda L m^2 \left\| \psi(\alpha^k) \right\| \right) \right] \left\| \psi(\alpha^k) \right\|.$$

Letting $\bar{\lambda}_k = \dfrac{1-\varepsilon}{Lm^2 \left\| \psi(\alpha^k) \right\|}$, we have

$$\left\| \psi(\alpha^k + \lambda p^k) \right\| \leq [1 - \varepsilon \lambda] \left\| \psi(\alpha^k) \right\| \quad (30)$$

for every $\lambda \in (0, \min\{1, \bar{\lambda}_k\})$. Then, (28), the definition of $\lambda_k$, implies that

$$\lambda_k \geq \min\{1, \xi \bar{\lambda}_k\}. \quad (31)$$

Since it follows from (30) that

$$\left\| \psi(\alpha^{k+1}) \right\| \leq [1 - \varepsilon \lambda_k] \left\| \psi(\alpha^k) \right\|, \quad (32)$$

$\{\left\| \psi(\alpha^k) \right\|\}$ is a monotonically decreasing sequence and so $\bar{\lambda}_k$ increases monotonically.

From (31), it follows that

$$1-\varepsilon\lambda_k \leq 1-\varepsilon\min\{1,\xi\overline{\lambda}_k\}\leq 1-\varepsilon\min\{1,\xi\overline{\lambda}_0\}.$$

Letting $q=1-\varepsilon\min\{1,\xi\overline{\lambda}_0\}$, then $q<1$ and by (32) we have

$$\|\psi(\alpha^k)\|\leq q^k\|\psi(\alpha^0)\|. \qquad (33)$$

Therefore, $\{\|\psi(\alpha^k)\|\}$ converges to zero with linear rate, which means that $\{\alpha^k\}$ converges to a solution of $\psi(\alpha)=0$ for any starting point $\alpha^0\in\Omega$. Since $\|\psi(\alpha^k)\|\to 0\,(k\to\infty)$, we have $\overline{\lambda}_k\to\infty\,(k\to\infty)$ and so there exists such $k_0$ that $\lambda_k=1$ for every $k\geq k_0$. Thus, in this case, (27) becomes the Newton method and the rate of convergence in the neighborhood of the solution is quadratic by the Lipschitz property (25). □

Let us consider sufficient condition under which the assumptions of Theorem 3 are satisfied. If

$$|h_i(x(\alpha))-h_i(x(\alpha'))|\leq L_i\|x(\alpha)-x(\alpha')\| \text{ and } \|x(\alpha)-x(\alpha')\|\leq M\|\alpha-\alpha'\|,$$

then $\|\psi'(\alpha)-\psi'(\alpha')\|\leq L\|\alpha-\alpha'\|$, because

$$\psi'(\alpha)-\psi'(\alpha')=\begin{pmatrix} h_1(x(\alpha))-h_1(x(\alpha')) & \cdots & 0 & 0 & \cdots & 0 \\ \vdots & \ddots & \vdots & \vdots & \ddots & \vdots \\ 0 & \cdots & h_N(x(\alpha))-h_N(x(\alpha')) & 0 & \cdots & 0 \\ 0 & \cdots & 0 & h_1(x(\alpha))-h_1(x(\alpha')) & \cdots & 0 \\ \vdots & \ddots & \vdots & \vdots & \ddots & \vdots \\ 0 & \cdots & 0 & 0 & \cdots & h_N(x(\alpha))-h_N(x(\alpha')) \end{pmatrix}.$$

Thus, the Lipschitz continuity of $x(\alpha)$ and $h(x)$ implies the Lipschitz continuity of $\psi'(\alpha)$. It is easy to see

$$[\psi'(\alpha)]^{-1}=\begin{pmatrix} \frac{1}{h_1(x(\alpha))} & \cdots & 0 & 0 & \cdots & 0 \\ \vdots & \ddots & \vdots & \vdots & \ddots & \vdots \\ 0 & \cdots & \frac{1}{h_N(x(\alpha))} & 0 & \cdots & 0 \\ 0 & \cdots & 0 & \frac{1}{h_1(x(\alpha))} & \cdots & 0 \\ \vdots & \ddots & \vdots & \vdots & \ddots & \vdots \\ 0 & \cdots & 0 & 0 & \cdots & \frac{1}{h_N(x(\alpha))} \end{pmatrix}$$

Thus, taking account of $\dfrac{1}{h_i(x(\alpha))}\leq u_i^u, i=1,...,N$, we have

$$\exists m,\ \forall\alpha\in\Omega,\ \left\|[\psi'(\alpha)]^{-1}\right\|\leq m,$$

that is, (26) is satisfied.

The modified Newton method defined by

$$\alpha^{k+1} = \alpha^k - \lambda_k \left[ \psi'(\alpha^k) \right]^{-1} \psi(\alpha^k)$$

can be rewritten component-wise as following.

$$\beta_i^{k+1} = \beta_i^k - \frac{\lambda_k}{h_i(x^k)} \left[ \beta_i^k h_i(x^k) - f_i(x^k) \right] = (1-\lambda_k)\beta_i^k + \lambda_k \frac{f_i(x^k)}{h_i(x^k)}, \ i=1,...,N \quad (34)$$

$$u_i^{k+1} = u_i^k - \frac{\lambda_k}{h_i(x^k)} \left[ u_i^k h_i(x^k) - 1 \right] = (1-\lambda_k)u_i^k + \lambda_k \frac{1}{h_i(x^k)}, \ i=1,...,N \quad (35)$$

From the above consideration, we construct an algorithm to find global solution of the problem (1) as following.

**[Algorithm MN]**

**Step 0.** Choose $\xi \in (0,1)$, $\varepsilon \in (0,1)$ and $y^0 \in X$. Let $\beta^0 = \frac{f(y^0)}{h(y^0)}, u^0 = \frac{1}{h(y^0)}, k=0$.

**Step 1.** Find a solution $x^k = x(\alpha^k)$ of the problem

$$\min \sum_{i=1}^N u_i^k \left( f_i(x) - \beta_i^k h_i(x) \right), \quad (36)$$

$$\text{subject to } x \in X$$

for $\alpha^k = (\beta^k, u^k)$.

**Step 2.** If $\psi(\alpha^k) = 0$, then $x^k$ is a global solution and so stop the algorithm. Otherwise, let $i_k$ denote the smallest integer among $i \in \{0,1,2,...\}$ satisfying

$$\left\| \psi(\alpha^k + \xi^i p^k) \right\| \leq (1-\varepsilon \xi^i) \left\| \psi(\alpha^k) \right\|$$

and let

$$\lambda_k = \xi^{i_k}, \ \alpha^{k+1} = \alpha^k + \lambda_k p^k, \ p^k = -\left[\psi'(\alpha)\right]^{-1} \psi(\alpha^k).$$

**Step 3.** Let $k = k+1$ and go to step 1. □

The step 2 of the algorithm MN is just the modified Newton method (27), (28). Therefore, the detailed update formula for $\alpha^k = (\beta^k, u^k)$ are (34) and (35).

If the stepsize $\lambda_k = \xi^{i_k}$ in the step 2 is replaced by $\lambda_k \equiv 1$, the step 2 is just the Newton method (24) and so we denote the algorithm MN by algorithm N in this case.

**Remark.** Suppose that we consider the following maximization problem instead of the problem (1).

$$\max F(x) = \sum_{i=1}^N F_i(x), \quad (37)$$

$$\text{subject to } g_i(x) \leq 0, i=1,\cdots,m, \ x \in R^n$$

where $F_i(x) = \frac{f_i(x)}{h_i(x)}$, $i = 1,...,N$, and $f_i(x)$, $g_i(x)$, $-h_i(x)$ are concave functions, respectively.

It is assumed that $f_i(x) \geq 0$ and $h_i(x) > 0$ in the feasible set $X = \{x \in R^n \mid g_i(x) \leq 0, i = 1,...,m\}$. Then, all of the results obtained above for the minimization problem (1) hold true for (37) except for that the problem (10) and (36) are replaced with the maximization problems.

## 4. Numerical experiments

Our algorithms have been implemented by 1.5GHZ 20GB 256MB PC NEC in the Matlab7.8 environment using the optimization toolbox. We made use of $|\Psi(\alpha^k)| \leq 1e-6$ as the stopping criterion in experiments below.

**Problem 1.**

$$\max \quad \frac{x_1}{x_1^2 + x_2^2 + 1} + \frac{x_2}{x_1 + x_2 + 1}$$

$$\text{subject to} \quad x_1 + x_2 \leq 1, \ x_1 \geq 0, \ x_2 \geq 0$$

The optimal value of the problem 1 is $f^* = 0.5958$. Our experiments had been carried out by both the algorithm N and algorithm MN: (i) simple iterative method (23), i.e. Newton method (24); (ii) modified Newton method (27) and (28), i.e. (34) and (35).

1) Starting point $x^0 = (0,0)^T$.
    (i) amount of iterations- 7,
    (ii) amount of outer iterations-7, amount of total iterations- 8.
2) Starting point – random point such that $x_1 + x_2 = 1, x_1 \geq 0, x_2 \geq 0$.
    The average performance of 100 runs are as follows:
    (i) amount of iterations- 7;
    (ii) amount of outer iterations- 7, amount of total iterations- 9.

In all experiments, we obtained the global solution with the accuracy of four-digits after 5 iterations. The algorithm MN required solving two subproblems in the first iteration and its stepsize had been always equal to 1 since the second iteration, that is, the algorithm MN turned to the algorithm N from the second iteration.

**Problem 2.**

$$\max \quad \frac{x_1}{x_1^2 + 1} + \frac{x_2}{x_2^2 + 1}$$

$$\text{subject to} \quad x_1 + x_2 \leq 1, \ x_1 \geq 0, \ x_2 \geq 0$$

The optimal solution and optimal value are $x^* = \left(\frac{1}{2}, \frac{1}{2}\right)^T$ and $f^* = \frac{4}{5}$, respectively.

1) Starting point $x^0 = (0,0)^T$
    (i) amount of iterations- 2,
    (ii) amount of outer iterations- 9,
        amount of total iterations- 16.
2) Starting point- random point such that $x_1 + x_2 = 1, x_1 \geq 0, x_2 \geq 0$
    The performance of 100 runs are as follows:

(i) success-0, failure-100 (local solutions);
(ii) average amount of iterations-4,
average amount of total iterations- 8.

**Problem 3.** We carried out numerical experiments for randomly generated problems as follows. The numerator and denominator of each term in objective function are

$$f_i(x) = \frac{1}{2}x^T A_{0i} x + q_{0i}^T x, \ i = 1,...,N$$

and

$$h_i(x) = c_i^T x, \ i = 1,...,N,$$

respectively. The feasible set is given by

$$X = \{x \in R^n \mid Ax \leq b, 1 \leq x_i \leq 5, i = 1,...,n\},$$

where $A_{0i} = U_i D_{0i} U_i^T$, $U_i = Q_1 Q_2 Q_3$, $i = 1,...,N$,

$$Q_j = I - 2\frac{\omega_j \omega_j^T}{\|\omega_j\|^2}, \ j = 1, 2, 3$$

and

$$\omega_1 = -i + rand(n,1), \ \omega_2 = -2i + 2 \cdot rand(n,1), \ \omega_3 = -3i + 3 \cdot rand(n,1),$$

$$c_i = i - i \cdot rand(n,1), q_{0i} = i + i \cdot rand(n,1),$$

$$A = -1 + 2 \cdot rand(5,n), \ b = 2 + 3 \cdot rand(5,1)$$

Starting with randomly generated point in $[1, 5]^n$, we carried out 200 runs of the algorithm N and the algorithm MN for fixed $n$ and $N$. Both of the algorithm N and the algorithm MN were successful for all generated problems. The average performance is shown in the table below.

(i) algorithm N, i.e.(24)

| n | N=5 | | N=10 | | N=50 | | N=100 | | N=200 | |
|---|---|---|---|---|---|---|---|---|---|---|
| | iter | time(s) | iter | time(s) | iter | time(s) | iter | time(s) | iter | time(s) |
| 10 | 6 | 0.2731 | 6 | 0.2753 | 6 | 0.2776 | 5 | 0.2844 | 5 | 0.2971 |
| 50 | 5 | 1.3288 | 5 | 1.3456 | 5 | 1.3764 | 5 | 1.7970 | 5 | 1.9313 |
| 100 | 5 | 4.4481 | 5 | 4.6139 | 5 | 4.6573 | 5 | 5.3741 | 5 | 5.6703 |
| 200 | 5 | 27.8647 | 5 | 27.9325 | 5 | 28.2704 | 6 | 29.1223 | 6 | 29.2663 |

(ii) algorithm MN, i.e. (34),(35)

| n | N=5 | | N=10 | | N=50 | | N=100 | | N=200 | | |
|---|---|---|---|---|---|---|---|---|---|---|---|
| | out | total | out | total | out | total | out | total | out | total | time(s) |
| 10 | 5 | 6 | 5 | 6 | 5 | 6 | 4 | 5 | 4 | 5 | 0.2984 |
| 50 | 5 | 6 | 4 | 5 | 4 | 5 | 4 | 5 | 4 | 5 | 1.1442 |

| 100 | 5 | 6 | 4 | 5 | 4 | 5 | 4 | 5 | 4 | 5 | 5.3757 |
| 200 | 5 | 6 | 4 | 5 | 5 | 6 | 4 | 5 | 4 | 5 | 29.4302 |

As shown in the above table, the amount of iterations of the algorithm N and the algorithm MN was independent of the amount of variables and fractional functions, and the run-time was proportional to the amount of variables. The experiments demonstrated that the algorithm N was successful for all problems but few problems , while the algorithm MN always gave global solution with any starting point.

## 6. Discussions

In this paper, we have presented a global optimization algorithm for the sum-of-ratios problem, which is non-convex fractional programming problem arising in computer graphics and finite element method. Our approach is based on transforming the sum-of-ratios problem into parametric convex programming problem and applying Newton-like algorithm to update parameters. The algorithm has global linear and local superlinear/quadratic rate of convergence. We have demonstrated that proposed algorithm always finds the global optimum with a fast convergence rate in practice by numerical experiments. In fact, the experiments showed that the optimal point is usually obtained within the first few iterations of the algorithm. So, we expect that our approach will solve successfully the several problems in multiview geometry formulated in [3].

## References

1. Carvalho, P. C. P., Figueiredo, L. H., Gomes, J., Velho, L., Mathematical Optimization in Computer Graphics and Vision, Elsevier, Morgan Kaufmann Publishers (2008).
2. Munson, T., Mesh shape-quality optimization using the inverse mean-ratio metric, Math. Program., 110, 561–590 (2007).
3. Kahl, F., Agarwal1, S., Chandraker, M. K., Kriegman, D., Belongie, S., Practical Global Optimization for Multiview Geometry, Int. J. Comput. Vis. 79,  271-284 (2008).
4. Freund, R.W., Jarre, F., Solving the sum-of-ratios problem by an interior-point method, J. Glob. Opt., 19, 83-102( 2001).
5. Tawarmalani, M., Sahinidis, N.V., Semidefinite relaxations of fractional programs via novel convexification techniques, J. Glob. Opt., 20, 137-158 (2001).
6. Benson, H.P., Using concave envelopes to globally solve the nonlinear sum-of-ratios problem, J. Glob. Opt., 22, 343-364 (2002).
7. Schaible, S., Shi, J., Fractional programming: the sum-of- ratios case, Opt. Meth. Soft., 18, 219-229 (2003).
8. Hartley, R., Sturm, P.: Triangulation, CVIU 68, 146–157 (1997).
9. Stew´enius, H., Schaffalitzky, F., Nist´er, D.: How hard is three-view triangulation really? In: Int. Conf. Computer Vision. 686–693 (2005).
10. Kahl, F., Henrion, D.: Globally optimal estimates for geometric reconstruction problems. In: Int. Conf. Computer Vision, Beijing, China, 978–985 (2005).
11. Ke, Q., Kanade, T.: Robust $L_1$- norm factorization in the presence of outliers and missing data by alternative convex programming. In: CVPR. 739–746 (2005).
12. Kahl, F.: Multiple view geometry and the $L_1$-norm. In: Int. Conf. Computer Vision, Beijing, China, 1002–1009 (2005).